# SCHOOL BUS ROUTING BY MAXIMIZING TRIP COMPATIBILITY

**Ali Shafahi**
PhD Candidate
Department of Civil and Environmental Engineering
University of Maryland
College Park, MD 20742
Email: ashafahi@umd.edu

**Zhongxiang Wang**
Graduate Student
Department of Civil and Environmental Engineering
University of Maryland
College Park, MD 20742
Email: zxwang25@umd.edu

**Ali Haghani**
Professor
Department of Civil and Environmental Engineering
1173 Glenn L. Martin Hall
University of Maryland
College Park, MD 20742
Phone: (301) 405-1963, Fax: (301) 405-2585
Email: haghani@umd.edu

Word count: 5,690 words text + (2 tables+ 5 figures) x 250 words (each) = 7,440 words

Submission Date July 29, 2016

**ABSTRACT**
School bus planning is usually divided into routing and scheduling due to the complexity of solving them concurrently. However, the separation between these two steps may lead to worse solutions with higher overall costs than that from solving them together. When finding the minimal number of trips in the routing problem, neglecting the importance of trip compatibility may increase the number of buses actually needed in the scheduling problem. This paper proposes a new formulation for the multi-school homogeneous fleet routing problem that maximizes trip compatibility while minimizing total travel time. This incorporates the trip compatibility for the scheduling problem in the routing problem. Since the problem is inherently just a routing problem, finding a good solution is not cumbersome. To compare the performance of the model with traditional routing problems, we generate eight mid-size data sets. Through importing the generated trips of the routing problems into the bus scheduling (blocking) problem, it is shown that the proposed model uses up to 13% fewer buses than the common traditional routing models.

*Keywords*: school bus routing, trip compatibility, school bus scheduling, bus blocking

# INTRODUCTION

School bus routing and scheduling problem (SBRSP) is traditionally broken into two pieces due to the extra complexity added when they are solved together. The first piece is school bus routing (SBRP) or trip building and the second piece is scheduling or blocking (SBSP). SBRP's goal is to build a *trip* that consists of several stops. PM trips usually start from schools and drop off students at different bus stop locations. AM trips are the reverse. SBSP or blocking, on the other hand, ties the trips together and assigns the grouped trips to *routes* (or *blocks*). In the blocking problem, it is possible to have buses that serve trips from different schools. If the school is a public school, the county's department of pupil transportation is the authority in charge of the operation and design of the school bus routes. From a monetary point of view, they are interested in having safe and reliable transportation services while keeping the cost for these services as low as possible. Since the major driving factor of cost is the number of buses (blocks), it is to their benefit to have buses with mixed trips.

In a private school setting, if the school owns and operates school buses, all trips of the blocks are for that school. In such setting, solving the routing problem and the scheduling problem independently does not cause any loss in terms of funds. This separation of the problem is thus only beneficiary as it improves the running time required for the search of the optimal solution. In the public school setting, however, the separation of these problems could result in financial losses and overall requirement of more school buses. The number of buses required for the county's operation is an output of the blocking (scheduling) problem. In the blocking problem, the compatibility of trips is the main influencing factor that helps reduce the overall number of buses. We say trip $b$ is compatible with trip $a$ if after trip $a$ is served, the bus has enough time to drive to trip $b$'s initial stop. If trips $a$ and $b$ are compatible, we can put them in the same block. However, if they are not compatible, they are assigned to separate blocks and this potentially increases the required number of blocks. The compatibility of the trips generated from the routing problem and the trips themselves are essential inputs of the blocking problem. Therefore, it is important that the school bus routing problems somehow generate more compatible trips.

The majority of the school bus routing problems, either have minimizing the number of total trips or minimizing the total travel time as their objective. These objectives are somewhat treating the routing problem and blocking problem as independent problems. Consider the example depicted in FIGURE 1. In this example, the three stops $s1, s2,$ and $s3$ are school bus stops for school 1. The population of school 1 exceeds the capacity of one bus and therefore we need at least two trips for that school. In addition, school 2 ends $2t$ times after school 1. Without considering the blocking, the minimum number of trips needed for school 1 is 2. The minimum total duration, $4t$, is achieved when trip 1 is: $School1 \rightarrow s1 \rightarrow s2$ and trip 2 is: $School1 \rightarrow s3$.. The soonest that any one of the trips can reach school 2 is at time $3t$. This is past the end time of school 2 and therefore, neither trip 1 nor trip 2 can be compatible with any of the trips of school 2.

Now, consider the case where trip 1 is: $school1 \rightarrow s1$ and trip 2 is: $school1 \rightarrow (s1) \rightarrow s2 \rightarrow s3$. The second trip does not stop at stop $s1$. It just passes through it to reach stop 2. In this case the total travel time is $5t$. However, trip 1 could reach school 2 at $2t$ and therefore it is compatible with the trips that are built for school 2.

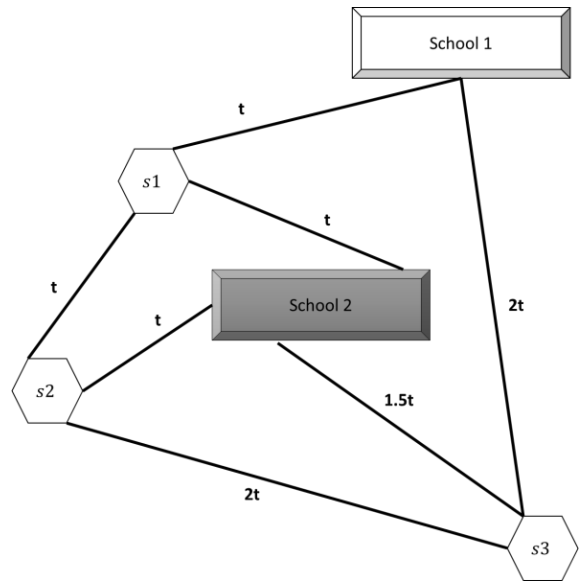

**FIGURE 1 Example that illustrates the importance of treating the school bus routing and scheduling problems as related problems**

Most existing research solve either SBRP or SBSP. They do not solve them together as it becomes computationally burdensome to find good solutions. However, the choice of a good objective for the school bus routing problem, that considers the SBSP (bus blocking), can potentially reduce cost and increase efficiency.

**Paper Structure**
In the next section, a review on school bus routing problem is presented. Then, the new mathematical model is explained. In the next section, computational tests are conducted to evaluate the model's performance. Finally, the performance of the proposed model, its applications, limitations, and future research directions are presented.

## LITERATURE REVIEW
The whole process of the school bus routing and scheduling problem involves five detailed steps: data preparation, bus stop selection, bus route generation, school bell time adjustment and route scheduling (1). The first two steps can be incorporated and solved by location-allocation-routing (LAR) strategy or allocation-routing-location (ARL) strategy. The difference and application of these two strategies can be referred to (2). The last two steps can be solved as a scheduling problem with time windows. School bus routing problem (SBRP) is always regarded as a variation of Vehicle Routing Problem (VRP). One difference between SBRP and classic VRP is that the travel time from the depot to the first pickup vertex (AM trips) is insignificant (2) as well as the travel time from the last drop-off vertex to the depot (PM trips). School bus service providers are mainly concerned about the ride time for students as opposed to total travel time for the buses. This makes SBRP become very similar to the Open Vehicle Routing Problem (OVRP). More specifically, many consider it to be a capacitated or distance (travel time) constrained OVRP (3). The main difference between OVRP and general VRP is that the former tries to find a set of Hamiltonian paths rather than Hamiltonian cycles as for the classical VRP (4). However, the minimum Hamiltonian path problem is still a NP-hard problem because it can be transferred into a minimum Hamiltonian cycle problem, which is a well-known NP-hard problem (5). Because this paper focuses on SBRP, which is similar to VRP, some VRP papers are also listed in the literature review.

**Characteristics of SBRP**

The school bus routing problem has been divided in various ways into different categories. Some divide the problem based on the number of schools for which the problem is solved. Another category is based on the area into which that the school(s) fall. The type of fleet is another mean for categorization.

In common public school settings, multi-school routing and scheduling is more realistic. However, quite a lot of papers focus on single-school problem due to its simplicity and similarity with classic VRP. When dealing with multi-school problems, many papers divide the problem into single-school problems by assuming each bus trip is exclusive for one school (6). Still, some papers consider mixed load service for multi-school problem. Braca et al. (7) solved the SBRP for multi-school bus routing in New York City with a Location Based Heuristic method, which forms routes by inserting the vertex with minimal insertion cost among all un-routed vertices and repeating this procedure by starting at random vertices. They then pick the best solution among all iterations.

Another classification is urban school and rural area school bus routing problem. In urban areas, where there tends to be more students in each stop, the bus capacity is usually the binding constraint (7) and stops may need to be served more than once. Therefore, maximum ride time constraint can be relaxed under certain conditions (8). While in rural areas, where each stop tends to have small number of students and the distance between adjacent stops are longer, the maximum ride time is usually a critical constraint. Moreover, vans or smaller vehicles rather than buses might be more economical, thus, vehicle type selection or mix fleet is preferred (9).

Homogeneous or heterogeneous fleet becomes important as it affects the bus capacity, the degree of crowding, the allowance of standing etc. National Association of State Directors of Pupil Transportation Services (10) regulates that a maximum of three young students (lower than the third grade) can sit in a typical 39-inch school bus seat.

(8) applied 3E (Efficiency, Effectiveness and Equity) criteria proposed by Savas (11) in SBRP. In order to balance the total ride time for each student, afternoon trips could be a replicated sequence (except the school) of the morning trip. Thus, it may be monetarily and temporally efficient to solve a PM (or AM) routing problem and then replicate it for the AM (or PM) trips. Balance of maximal load or maximal ride time is another equity concern, which has been considered as an objective in the formulation (8, 12).

**Objective**

The most common objectives of the SBRP are minimizing the total travel distance or travel time (8, 9, 13-24), minimizing the total number of trips (6, 7, 14,15, 25) or minimizing the total cost including bus purchase (fixed) cost and bus operation cost (3, 26, 27, 28). In the latter, the total cost is a combination of the first two objectives considering bus purchase cost is proportional to the number of trips (without scheduling, each trip is assumed to be served by one bus) and bus operation cost is equivalent to the total travel distance. (8) incorporated students' walking distance in their model. (6) and (25) both tried to minimize the number of buses while minimizing the maximum ride time for students. These two objectives are conflicting because less number of buses would require longer route lengths (6). Another factor in the objective function could be penalty. Penalties are generally added when a hard constraint is relaxed and we allow violation of constraints. This penalty could have various forms including time window violation penalty (19), and capacity or maximum ride time constraints violation penalty (20).

**Constraints**

Logistic (or degree) constraints are the basic VRP constraints, which formulate trips. They regulate the conservation of flow by assuring that when a vehicle enters a vertex (node), it also exits the same vertex. In addition, each trip starts from a depot and goes back to the depot, each vertex is served exactly once for single-visit cases. Multi-visit is allowed in some cases, but most papers only consider single-visit scenario. Although SBRP is an OVRP as discussed before, it is easier to generate round trips first and then discard the first (AM trips) or last (PM trips) arc of that trip. Capacity constraints are also quite common. All papers about SBRP and most papers about VRP have vehicle capacity constraints. Maximum ride time constraint is also considered, especially in rural areas. Russell and Morrel (13) used 45 minutes as the maximum ride time and Chen et al. (14) limited the ride time to be less than or equal to 75 minutes.

Sub-tour elimination constraints are another common set of constraints that prevent the formation of the illegal trips some of which are not connected to the depot (3). The most adopted way to enforce connectivity is to set an appropriate lower bound on the number of vehicles required to visit a subset of vertices (29). Another common method is using an artificial commodity flow and the introduction of flow variables (8). Bektaş and Elmastaş (3) adopted Miller, Tucker and Zemlin (30) subtour elimination constraint to solve SBRP.

Chain Barring constraints are only seen in multi-school scenarios. These constraints aim to eliminate the illegal trips that start from one school and end at another school (26). These constraints can be enforced implicitly by simply allowing trips assigned to a school to only stops that are assigned to the same school. Thus, chain barring constraints are considered implicitly in this paper. We should note that an assumption of this implied constraint is that we cannot have a school trip with mixed students from different schools. Yet, this is not a huge shortcoming as in reality mixed students are not popular among school agencies.

Time windows is widely incorporated in scheduling rather than route generation. (7) applied time window constraints as school buses must arrive at a school no earlier than 25 minutes before and no later than 5 minutes before the start of school. (28) proposed time window constraints, which prevent the sub-tour formation at the same time.

For solving routing and scheduling together, (31) modified several routing heuristic algorithms (saving algorithm, insertion algorithm, nearest neighbor algorithm, and the sweep algorithm) to solve the routing and scheduling problem. However, the optimality of these heuristics cannot be proven nor can they be evaluated due to the lack of optimal solutions from exact algorithms.

In this paper, we propose a new mathematical model with a new bi-objective objective function. The proposed model maximizes trip compatibility while minimizing total travel time. It does not solve the routing and scheduling problem simultaneously and hence finding good solutions for the problem is not computationally burdensome. However, it provides a routing solution (a set of trips) with high compatibility. As a result, multi-trip routes that can be served by one bus is more likely to be generated and number of buses that are needed for the entire school operations could be reduced.

## METHODOLOGY

We present a Mixed Integer Programming (MIP) model that aims to maximize the compatibility among trips with the hope of minimizing the overall number of buses as a result. The presented MIP is a novel model that solves the routing problem while considering the bus blocking problem. TABLE 1 summarizes the variables and parameters used in the mathematical model.

**TABLE 1 Notation summary of variables and parameters**

| | |
|---|---|
| Variables for school bus routing | |
| Variable name | Description |
| $s2t_{s,t}$ | Binary decision variable that equals 1 if stop $s$ is assigned to trip $t$ |
| $t2s_{t,k}$ | Binary decision variable that equals 1 if trip $t$ is assigned to school $k$ |
| $p4t_{s,t}$ | Portion of the capacity of the bus doing trip $t$ that is filled at stop $s$ |
| $x^t_{s1,s2}$ | Binary variable that equals 1 if in trip $t$ the bus goes directly from stop $s1$ to stop $s2$ |
| $tt_t$ | Travel time (duration) of trip $t$ |
| $start_t$ | The start time of trip $t$ |
| $end_t$ | The end time of trip $t$ |
| $b_{t1,t2}$ | Binary variable that equals 1 if trips $t2$ can be served after trip $t1$ (they are compatible) |
| $dd_{t1,t2}$ | The travel deadhead duration from the last stop of trip $t1$ to the first stop of trip $t2$ |
| $l_{s,t}$ | Binary variable that equals 1 if the last stop of trip $t$ is stop $s$ |
| $c^t_{s1,s2}$ | Variable used for subtour elimination – The units of "artificial commodity" that is shipped from stop $s1$ to $s2$ using trip $t$ |
| Variables for bus blocking | |
| $y_{t1,t2}$ | Binary variable that equals 1 if trip $t2$ follows trip $t1$ within the same block (tour). This variable is only created for compatible pair of trips. |
| $a_t$ | Binary variable that equals 1 if a trip belongs to a tour of size 1 (no trip goes before or after it) |
| $m_t$ | A binary variable that equals 1 if trip $t$ has a trip that is being served before it and also another trip that is being served afterwards |
| $f_t$ | A binary variable that is equal to 1 if trip $t$ is served as the first trip in a tour of size greater than 1 |
| Parameters for School bus routing and bus blocking | |
| Parameter name | Description |
| $Schools$ | Set of schools |
| $Trips\_School_k$ | Set of possible trips dedicated to school $k$ |
| $Trips$ | Set of all trips |
| $OTrips$ | The set of outputted trips from school bus routing problems |
| $Stops\_School_k$ | The set of stops in which students for school $k$ should go to / come from |
| $Stops$ | Set of all stops |
| $Bus\_Cap$ | The capacity of each bus |
| $Students_{s,k}$ | The number of students at stop $s$ for school $k$ |
| $O_k$ | Origin stop for school $k$ |
| $D_{s1,s2}$ | The duration to drive from stop $s1$ to $s2$ plus the dwell time required at the stops |
| $SchoolEnd_k$ | Time which school $k$ closes |
| $Buffer\_pickup$ | The allowable buffer for picking up students from the schools = 0 for cases solved due to the trips being PM trips with no flexibilities. |
| $M$ | A large value (big-M) |
| $C_B$ | Penalty for total number of trips = 1000 for the cases solved |
| $C_C$ | Coefficient of compatibility = 200 for the cases solved |
| $A$ | The number of allowable trips over the minimum required trips for each school |

As explained earlier, the main purpose of the objective function is to simultaneously maximize the number of compatible trips while minimizing the number of overall trips required and the travel time. These are done by assigning penalties (positive weights) to total number of trips being used and assigning benefits (negative weights) to compatible pair of trips.

$$\min z = \sum_{t \in Trips} tt_t + C_B \sum_{k \in Schools} \sum_{t \in Trips\_School_k} t2s_{t,k} + C_C \sum_{t1 \in Trips} \sum_{t2 \in Trips} b_{t1,t2} \quad (1)$$

The constraints listed below are in charge of building trips for each school. Each trip is built out of various stops.

$$s2t_{s,t} \leq t2s_{t,k} \quad \forall k \in Schools, s \in Stops\_School_k, t \in Trips\_School_k \quad (2)$$

$$\sum_{s \in School\_Stops_k} p4t_{s,t} \leq 1 \quad \forall k \in Schools, \forall t \in Trips\_School_k \quad (3)$$

$$\sum_{t \in Trips\_School_k} p4t_{s,t} = \frac{Students_{s,k}}{Bus\_Cap} \quad \forall k \in Schools, s \in School\_Stops_k \quad (4)$$

$$p4t_{s,t} \leq s2t_{s,t} \quad \forall k \in Schools, t \in Trips\_School_k, s \in School\_Stops_k \quad (5)$$

$$\sum_{s2 \in Stops\_School_{sch}|s2 \neq s} x^t_{s,s2} = s2t_{s,t} \quad \forall k \in Schools, \forall t \in Trips\_School_k, \forall s \in Stops\_School_k \quad (6)$$

$$\sum_{s \in Stops\_School_k} x^t_{O_k,s} = t2s_{t,k} \quad \forall k \in Schools, t \in Trips\_School_k \quad (7)$$

$$\sum_{j \in Stops\_School_k | j \neq s} x^t_{s,j} = \sum_{i \in Stops\_School_k | i \neq s} x^t_{i,s} \quad \forall k \in Schools, s \in Stops_{School_k}, \forall t \in Trips\_School_k \quad (8)$$

$$tt_t = \sum_{s1 \in Stops\_School_k} \sum_{s2 \in Stops\_School_k | s2 \neq O_k} x^t_{s1,s2} \times D_{s1,s2} \quad \forall k \in Schools, t \in Trips\_School_k \quad (9)$$

$$end_t = start_t + tt_t \quad \forall t \in Trips \quad (10)$$

$$SchoolEnd_k \leq start_t \leq SchoolEnd_k + Buffer\_pickup \quad \forall k \in Schools, t \in Trips\_School_k \quad (11)$$

$$end_{t1} + dd_{t1,t2} - M \times (1 - b_{t1,t2}) \leq start_{t2} \quad \forall t1, t2 \in Trips \quad (12)$$

$$dd_{t1,t2} \geq \frac{M}{2} \times (1 - t2s_{t1,k1}) + \frac{M}{2}(1 - t2s_{t2,k2}) + \sum_{s1 \in Stops\_School_{k1}} D_{s1,O_{k2}} \times l_{s1,t1} \quad \forall k1, k2 \in Schools, t1 \in Trips\_School_{k1}, t2 \in Trips\_School_{k2} \quad (13)$$

$$x^t_{s,O_k} = l_{s,t} \quad \forall k \in Schools, t \in Trips\_School_k, s \in Stops\_School_k \quad (14)$$

$$t2s_{t1,k} \geq t2s_{t2,k} \quad \forall k \in Schools, t1, t2 \in Trips\_School_k \mid t2 \leq t1 \quad (15)$$

$$\sum_{i \in Stops\_School_k | i \neq s} c^t_{i,s} - \sum_{i \in Stops\_School_k | i \neq s} c^t_{s,i} = s2t_{s,t} \quad \forall k \in Schools, t \in Trips\_School_k, s \in Stop\_School_k | s \neq O_k \quad (16)$$

$$c^t_{s1,s2} \leq M \times x^t_{s1,s2} \quad \forall k \in Schools, t \in Trips\_School_k, s1, s2 \in Stop\_School_k \quad (17)$$

$$\left\lceil \frac{\left(\sum_{s \in Stop_{School_k}} Students_{s,k}\right)}{Bus\_Cap} \right\rceil \leq \sum_{k \in Schools} \sum_{t \in Trips_{School_k}} t2s_{t,k} \qquad \forall k \in Schools \qquad (18)$$

$$\sum_{k \in Schools} \sum_{t \in Trips_{School_k}} t2s_{t,k} \leq A + \left\lceil \frac{\left(\sum_{s \in Stop_{School_k}} Students_{s,k}\right)}{Bus_{Cap}} \right\rceil \qquad \forall k \in Schools \qquad (19)$$

Constraints (2) prevent the assignment of stops to the trips that have not been assigned to their respective schools. These trips are the ones that are not in use. Constraints (3) are trip capacity constraints. Constraints (4) ensure that all students are being served. Constraints (5) disallow the assignment of students to a trip that does not pass the stop in which the students are from there. Constraints (6) assure that if a stop is assigned to a trip, that trip visits that node by traversing the arcs that lead to that node. Constraints (7) enforces each trip to start from the school for PM trips. Conservation of flow is expressed through Constraints (8).

Constraints (9) calculate the travel duration of the trips. The trips, start from the school and end at a node, which is the last stop (the stop before going back to the school). Constraints (10) calculate the end time of the trips using the start time and the travel time. The start time is constrained by the time at which the schools close in constraints (11). Constraints (12) are used for identifying the compatible trips. The deadhead between pairs of trips that is used for compatibility calculations is computed using constraints (13). Constraints (14) are used for identifying the last stops for the trips. This last stop is the last actual stop and is the stop right before the school goes back to its start point to complete its fictional closed loop trip. Constraints (15) are for eliminating symmetries. Note that if there is only one trip assigned to a school, the id for this trip could assume many values. These constraints prevent higher trip ids to occur prior to the lower ones and can speed up the search for good solutions. Constraints (16)-(17) are flow subtour elimination constraints. Finally, constraints (18)-(19) are limiting the number of trips assigned to a school to the minimum number needed based on population at each stop. As discussed before, in an urban school setting where bus capacity is the binding constraint, maximum ride time constraint can be relaxed for the sake of improving model efficiency.

To compare the trip outcomes from our model versus those from traditional models, we input the trips from all models into a bus blocking problem summarized below:

$$\sum_{j \in OTrips} y_{t,j} + \sum_{j \in OTrips} y_{j,t} + a_t \geq 1 \qquad \forall t \in OTrips \qquad (20)$$

$$\sum_{j \in OTrips} y_{t,j} + a_t \leq 1 \qquad \forall t \in OTrips \qquad (21)$$

$$\sum_{j \in OTrips} y_{j,t} + a_t \leq 1 \qquad \forall t \in OTrips \qquad (22)$$

$$\sum_{j \in OTrips} y_{t,j} + \sum_{j \in OTrips} y_{j,t} \leq 2 \qquad \forall t \in OTrips \qquad (23)$$

$$\sum_{j \in OTrips} y_{t,j} + \sum_{j \in OTrips} y_{j,t} + a_t = 1 + m_t \qquad \forall t \in OTrips \qquad (24)$$

Constraints (20)-(23) make sure that each trip that is selected and assigned to a school is either served with another trip as a pair or being served alone. Constraints (24) identify the middle trips of blocks. The objective of the blocking problem is minimizing the total number of blocks. This could be achieved by maximizing the number of middle trips, $m_t$, and minimizing the number of alone trips, $a_t$ (25). Alone trips are those trips that are the sole trips in a tour (block) of size one.

$$\min z = -\sum_{t \in OTrips} m_t + \sum_{t \in OTrips} a_t \tag{25}$$

It is important to note that the way the objective function is written in (1) may be putting too much importance on the total trip compatibility. Consider a scenario that trip a is compatible with both trips b and c but trips b and c are incompatible. In such a case, one bus can serve either trip b or trip c after trip a, but is no bus can serve them all. Total trip compatibility is equal to two while only one bus can be saved. This implies that trip compatibility is not equal to the number of buses that can be saved. One approach is to have a hybrid model of routing and blocking. However, as mentioned earlier, such a problem will be very complex due to the additional variables and constraints. Such a problem is almost impossible to solve for problems that are not toy problems using exact algorithms. Another approach is to look at compatibility of trips with schools. We have looked into this approach however, upon preliminary investigations, we have did not find any great gain over the proposed model.

## COMPUTATIONAL RESULT

In order to evaluate the model, eight set of randomly generated mid-size problems were generated and tested against typical routing models found from literatures, that are minimizing number of trips (which is same as minimizing number of buses without scheduling) and minimizing total travel time (or distance). The parameters of the distributions used for the randomly generating the cases are mostly based on real world cases. The detailed information for each scenario and the results from different models are shown in TABLE 2. Each scenario has many parameters. The total number of stops (# stops) is an indicator of the problem size, which is the number of vertices in the problem. The number of schools also affects the size of the problem since each school has a certain number of students that require bus transport. The total number of trips that each school needs is bound by the ratio of the total number of students over the capacity of each bus that is 48 for each school. The average number of students at each school that require bus transportation are also different for each scenario. Note that each school has a different population. The Maximum number of stops assigned to each school is another parameter that influences the problem's complexity. Finally, bus start range for each scenario is the time difference between the earliest and latest school dismissal time. The dismissal times for the schools follow a discrete uniform distribution with 15 minutes time intervals within that range.

For each scenario, four objectives are tested: 1) the proposed objective of maximizing trip compatibility while minimizing total travel time (**MaxCom+TT**); 2) maximizing trip compatibility (**MaxCom**) by deleting travel time from the proposed objective; 3) minimizing total number of trips (**MinN**); and 4) minimizing total travel time (**MinTT**). The problems are solved by commercial solver FICO Mosel XPRESS on five computers with same feature: Intel® Core™ i5-2400 CPU, 3.10 GHz with 4 GB RAM. Due to the slow rate of reduction in the gap, the solution processes for the test problems were terminated after a certain running times (from 30 minutes to 24 hours). The relative performance measure used for each of the models, is the total number of blocks that is the output of the blocking problem when the input trips are the solutions of the respective model.

Scenario 1 is for illustrating the impact of "additional allowed trips" (i.e. the A parameter in constraints (19)) on the complexity of the problem and solution time. To perform this sensitivity analysis, we run scenario 1 using four different values for the A parameter (A=0,1,2,3). Recall that if A=2, each school is allowed to use up to two more trips than the minimum number it requires. The minimum number of trips a school requires is easily calculated based on the school's population and the capacity of each bus. As it can be seen, the more additional trips are allowed, the more running time is required. However, better solutions may be found. For example, 30 buses (MaxCom+TT, 2 additional trips after 6 hrs.) is the current best solution in comparison to 3 additional trips with short running time (30 minutes) and all other solutions from 1 and 0 additional trips cases. While keeping in mind the recently mentioned note, for the sake of saving running time, this paper limits all the remaining 7 scenarios to have no additional trips (A=0). It should be noted that in real applications, additional trips should be allowed to try to find the best potential solution by simply allowing the models to run longer. In real applications, the running time for SBRSP is less important than finding the optimal solution. Even one bus saving could be significantly beneficial.

**TABLE 2 Summary of computational result**

| Scenario | | 1 | | | | | 2 | 3 | 4 | 5 | 6 | 7 | 8 |
|---|---|---|---|---|---|---|---|---|---|---|---|---|---|
| # of Stops | | 100 | | | | | 200 | 100 | 100 | 125 | 100 | 200 | 200 |
| # of Schools | | 5 | | | | | 20 | 20 | 20 | 25 | 20 | 20 | 20 |
| Avg. # of Student to School | | 91.4 | | | | | 89.6 | 120.7 | 182.8 | 90.4 | 91.6 | 89.5 | 91.1 |
| Max # of Stops to School | | 13 | | | | | 16 | 13 | 13 | 13 | 13 | 16 | 14 |
| Bus Service Start Time Range | | 0-30 | | | | | 0-30 | 0-30 | 0-30 | 0-30 | 0-90 | 0-90 | 0-16 |
| Additional Trips | 3 | 2 | 2 (6 hrs.) | 1 | 0 | 0 | 0 | 0 | 0 | 0 | 0 | 0 | 0 |
| Max Com +TT | RT | 30 | 30 | 180 | 30 | 30 | 30 | 60 | 600 | 360 | 30 | 60 | 30 |
| | Gap (%) | 21.1 | 22.9 | 6.7 | 8.4 | 5.0 | 6.2 | 5.3 | 12.0 | 2.3 | 2.4 | 7.2 | 1.9 |
| | NOT | 50 | 52 | 42 | 42 | 41 | 47 | 54 | 75 | 46 | 41 | 47 | 45 |
| | NOB | **34** | **34** | **30** | **30** | **31** | **32** | **34** | **45** | **35** | **16** | **24** | **45** |
| Max Com | RT | 30 | 30 | 180 | 30 | 30 | 30 | 60 | 600 | 360 | 30 | 60 | 30 |
| | Gap (%) | 10.7 | 25.6 | 7.3 | 16.1 | 4.9 | 5.6 | 4.9 | 10.8 | 1.5 | 2.7 | 10.5 | 1.3 |
| | NOT | 43 | 54 | 42 | 47 | 41 | 47 | 54 | 75 | 46 | 41 | 47 | 45 |
| | NOB | **31** | **34** | **31** | **32** | **30** | **30** | **36** | **45** | **36** | **16** | **24** | **44** |
| Min N | RT | 0.38 | 0.85 | 0.05 | 0.14 | 0.10 | 0.46 | 0.04 | 0.04 | 0.02 | 0.20 | 0.35 | 0.18 |
| | Gap (%) | 0 | 0 | 0 | 0 | 0 | 0 | 0 | 0 | 0 | 0 | 0 | 0 |
| | NOT | 41 | 41 | 41 | 41 | 41 | 47 | 54 | 75 | 46 | 41 | 47 | 45 |
| | NOB | **36** | **35** | **35** | **35** | **38** | **38** | **45** | **60** | **41** | **21** | **26** | **45** |
| Min TT | RT | 30 | 30 | 180 | 30 | 30 | 30 | 60 | 600 | 360 | 30 | 60 | 30 |
| | Gap (%) | 1.1 | 1.1 | 0.8 | 1.1 | 0.9 | 3.4 | 6.2 | 20.6 | 0.2 | 1.4 | 7.0 | 15.0 |
| | NOT | 61 | 58 | 59 | 58 | 41 | 47 | 54 | 75 | 46 | 41 | 47 | 45 |
| | NOB | **40** | **40** | **39** | **39** | **32** | **32** | **38** | **52** | **35** | **17** | **24** | **45** |

RT: Running time (minute); Gap: unit in % NOT; Number of trips (trip); NOB: Number of buses (bus)

In Scenario 2, where the number of stops is increased, the model presented in this paper can save 2 buses in comparison to the solution from MinTT and 8 buses can be saved from MinN (FIGURE 2). Considering the higher gap for proposed model under same running time (30 minutes), the result has the potential to be further improved. Scenario 3 is similar with scenario 1, except increasing average number of students for each school from 91.4 students to 120.7 students.

Based on the result from MaxCom+TT, 2 and 11 buses can be saved compared to that from MinTT and MinN, respectively. Scenario 4 is a further expansion of Scenario 1 and Scenario 3, which has an average of 182.8 students for each school. In this relatively large case, the saving by MaxCom+TT and MaxCom is more significant, 7 buses and 15 bused are saved. In Scenario 5, both number of stops and number of school are increased. In this scenario, MaxCom+TT and MinTT both need 35 buses, but MaxCom needs 36 buses. The reason might be the bigger gap for MaxCom. As it can be seen, in general, MaxCom+TT uses less buses compared to traditional models.

     Scenario 6 is similar with Scenario 1 except it expands the bus service start time range from 30 minutes to 90 minutes. The purpose of this change is to see how the model performs in comparison to traditional models when many of the trips are going to be easily compatible anyway. Because trips can be easily compatible with each other after increasing start time range, bus saving might be less than smaller time range scenarios. The result in FIGURE 2 proves this hypothesis, only one bus saved in Scenario 6. Scenario 7 is similar to Scenario 1 with a change in the bus service start time range. The range is increased to 90 minutes for Scenario 7. It has similar results to Scenario 1, no additional improvement made by proposed models. Scenario 8 is another experiment, which decreases the bus service start time range to 16 minutes from Scenario 2. Only one bus is saved by MaxCom. In these conditions, the bus service start time range is too small for trips to be compatible. MaxCom would take whatever it need to make it compatible, which leads to extremely long travel time trips (FIGURE 3h). The cross examination reveals the application and limits of proposed objective. That is the proposed objective is not extremely beneficial in comparison to the traditional models for scenarios where the bus service time range falls in a range that is either too big for trips being compatible anyway or too small that trips cannot be compatible. But for the range in between, like school bus service provided by a county's Department of Pupil Transportation, the proposed objective, MaxCom+TT greatly outperforms the traditional objectives.

     Another concern in SBRSP is the maximum ride time for each trip. Since the MaxCom model does not have a maximum ride time constraint nor any incentive in minimizing the total travel time, MaxCom only focuses on maximizing the trip compatibility even at the expense of long travel time trips (that could go up to 121 minutes in Scenario 4, which is too long to be a school bus trip). From this perspective, MaxCom+TT makes more sense. FIGURE 3 is the travel time distribution where frequency is calculated in 5 minutes interval and marked at the beginning of each interval. For instance, frequency of interval 0-5 minutes is plotted at travel time (horizontal axis) equals to 0 minutes. The results in FIGURE 3 show that MaxCom+TT tends to have short travel time trips (mostly less than 40 minutes) than MaxCom. FIGURE 4 is the total travel time for different objectives. It shows that in most scenarios, the results from MinTT has the minimal total travel time, and that MaxCom and MinN have long travel times. MaxCom+TT has much less total travel time than those from MaxCom and MinN but the total travel times are slightly higher than that from MinTT.

     In general, MinN has the worst result. It is because the model would stop right after finding a solution that uses the minimal number of trips. Since the objective cannot be improved furthermore (i.e. it is optimal), the model simply ignores decreasing travel time or increasing trip compatibility, which makes it harder to group them in a block. MamCom+TT tends to have more high-compatible (2- or 3-trips) routes (FIGURE 2). The more high-compatible routes exist, the less buses are needed. Overall, MaxCom+TT generates the best result, which use the least buses and its trips have relatively short travel times.

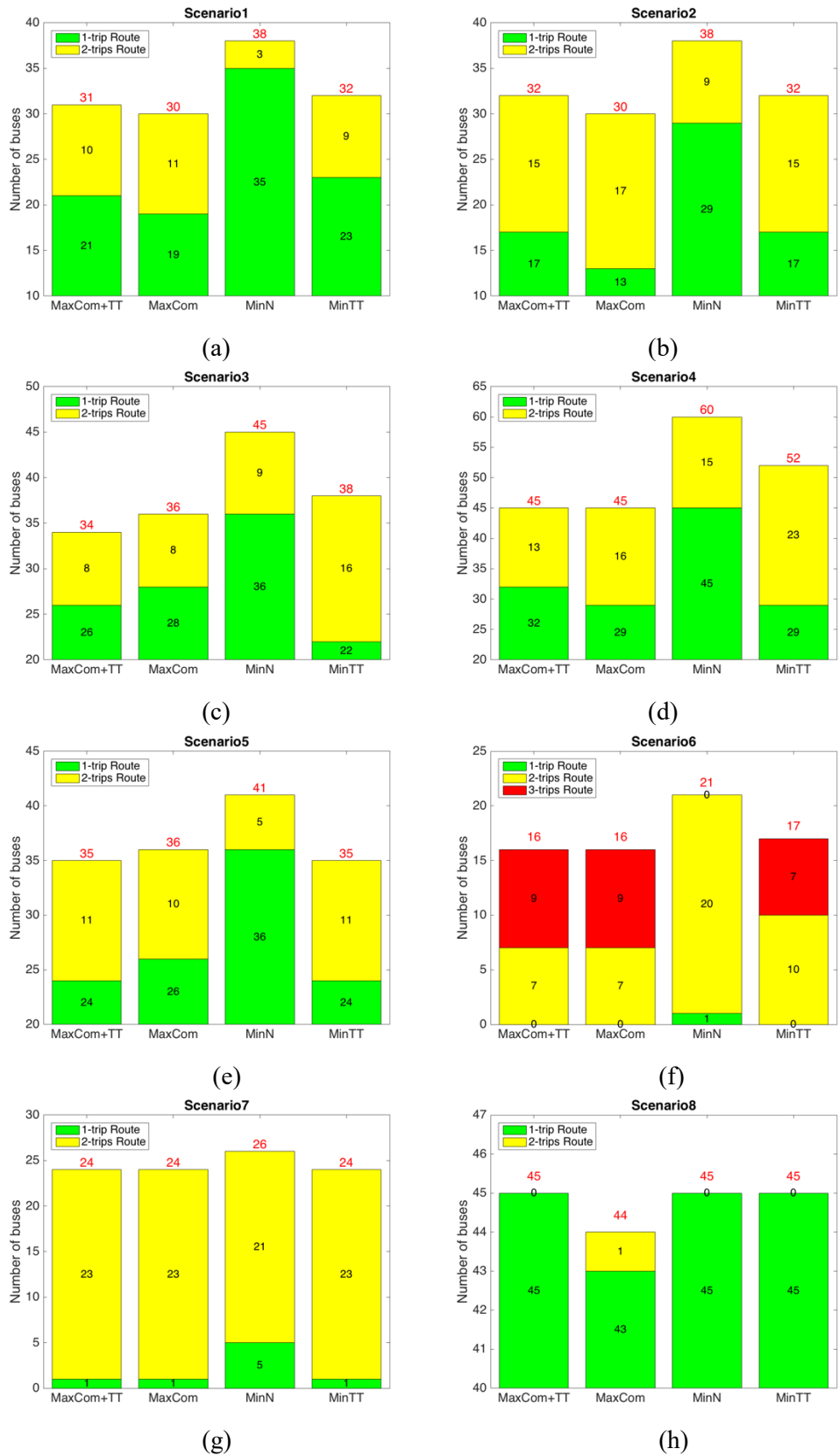

**FIGURE 2 Summary of school bus blocking result**

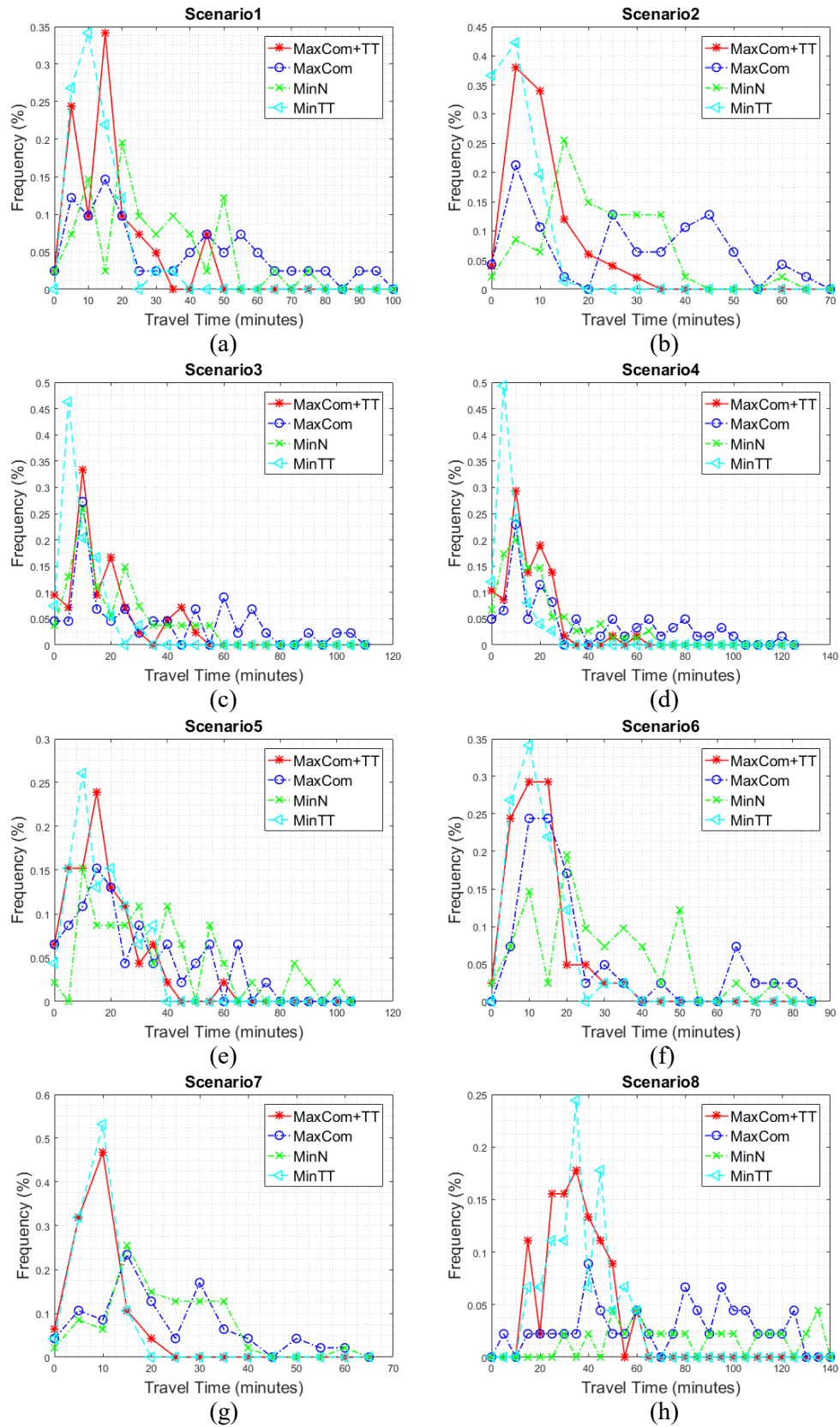

**FIGURE 3 Travel time distribution**

Shafahi, Wang, Haghani                                                                                                           15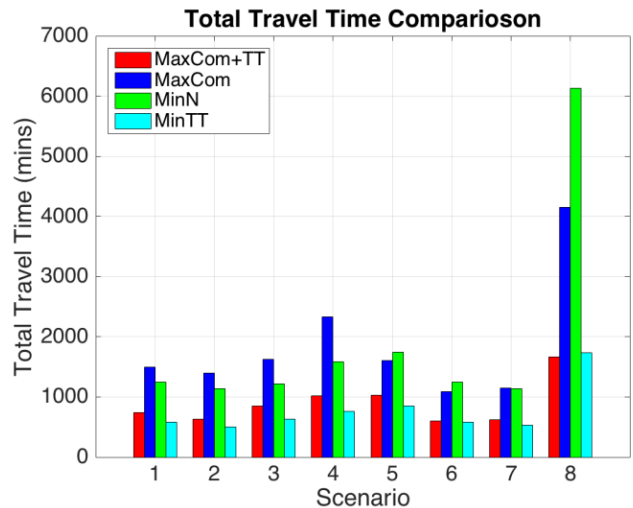

**FIGURE 4 Total travel time comparison**

It is worth noting that this gain in savings in the number of buses needed for serving the trips does not come free. FIGURE 5 shows the tradeoff between the reduction in the number of buses and the increase in travel time per bus. The comparison is made between the best results from the proposed models (MaxCom & MaxCom+TT) and the best results from traditional models (MinTT & MinN). It can be seen that up to 13.5% buses can be saved at the expense of 5.8 additional travel time minutes per bus (Scenario 4, 7 buses are saved). From a financial point of view, the savings gained by needing fewer buses could easily justify the additional travel times.

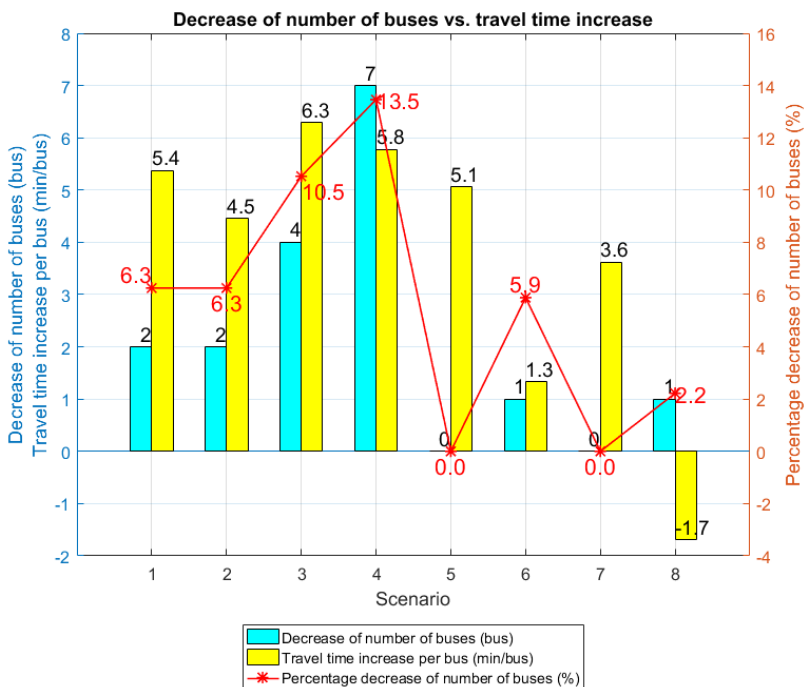

Note: -1.7 for Scenario 8 means travel time per bus of the solution from new objectives is 1.7 minutes less than that from traditional objectives.

**FIGURE 5 Tradeoff between decrease of number of buses and travel time increase from best new models and best traditional models**

## SUMMARY AND CONCLUSION

In this paper, we propose a new mathematical model that optimizes a new objective to solve school bus routing problem. This new objective has a component for maximizing the trip compatibility that can potentially decrease the total number of buses required for the operation. Eight mid-size problems are solved to test the performance, and illustrate the applications and limitations of the proposed method. It is shown that the proposed model outperforms traditional models by requiring fewer buses. The model's significance is greater when the school dismissal times are within 30 minutes of each other. We recommend using this model for such cases.

This research opens the venue for many future studies. One of them is sensitivity analysis for parameters in the proposed objective. Another one is the double counting of compatibility problem that was mentioned at the end of the methodology section. Larger size problem should also be applied to test the performance and efficiency of the proposed model. A heuristic algorithm should also be developed to solve the model for larger problems when needed.

Shafahi, Wang, Haghani 17